\documentclass[3p]{elsarticle}
\usepackage{amssymb, color}
\usepackage{amsmath}
 \usepackage{pifont}
 \usepackage[utf8]{inputenc}

\makeatletter
\def\LaTeX{\leavevmode L\raise.42ex
    \hbox{\kern-.3em\size{\sf@size}{0pt}\selectfont A}\kern-.15em\TeX}
\makeatother

\newcommand{\BibTeX}{{\rm B\kern-.05em{\sc
          i\kern-.025emb}\kern-.08em\TeX}}

\makeatletter
\def\@currentlabel{2.1}\label{e:dispaa}
\def\@currentlabel{2.21}\label{e:dispau}
\def\@currentlabel{2.22}\label{e:dispav}
\def\@currentlabel{2.23}\label{e:dispaw}
\def\@currentlabel{2.24}\label{e:dispax}
\def\theequation{\thesection.\@arabic\c@equation}
\makeatother

\renewcommand{\theequation}{\arabic{section}.\arabic{equation}}

\newcommand{\R}{\mathbb R}

\newcommand{\N}{\mathbb N}
\newcommand{\e}{\epsilon}

\def \D{\Delta}
\def \O{\Omega}
\newtheorem{thm}{Theorem} [section]
\newtheorem{lem}{Lemma} [section]
\newtheorem{prop}{Proposition} [section]

\newtheorem{rmk}{Remark}[section]
\newtheorem{taggedtheoremx}{Theorem}
\newenvironment{taggedtheorem}[1]
 {\renewcommand\thetaggedtheoremx{#1}\taggedtheoremx}
 {\endtaggedtheoremx}

\renewcommand{\theequation}{\thesection.\arabic{equation}}
\renewcommand{\thesection}{\arabic{section}}
\renewcommand{\theequation}{\thesection.\arabic{equation}}
\let\ssection=\section\renewcommand{\section}{\setcounter{equation}{0}\ssection}

\def \p{\partial}

\begin{document}
\begin{frontmatter}

\title{Explicit estimates for solutions to higher order elliptic PDEs via Morse index}

\author[hh]{Abdellaziz Harrabi }
\ead{abdellaziz.harrabi@yahoo.fr}
\author[fd]{Foued Mtiri}
\ead{mtirifoued@yahoo.fr}
\author[ah]{Dong Ye\corref{cor1}}
\ead{dong.ye@univ-lorraine.fr}
\address[hh]{Institut Sup\'erieur des Math\'ematiques appliqu\'ees et de l'Informatique. Universit\'e de Kairouan, Tunisia}
\address[fd]{IECL, UMR 7502, Universit\'{e} de Lorraine, France}
\address[ah]{IECL, UMR 7502, Universit\'{e} de Lorraine, France}
\begin{abstract}

In this paper, we establish $L^{\infty}$ and $L^{p}$ estimates for
solutions of some polyharmonic elliptic equations via the Morse index. As far as we know, it seems to be the first time 
that such explicit estimates are obtained for polyharmonic problems.
\end{abstract}
\end{frontmatter}

 \section{Introduction}
Consider the following polyharmonic equations $(P_k):$\,$(-\Delta)^k u=f(x,u)$ in $\O$
with the Dirichlet boundary conditions
\begin{align}
\label{D}
u = \frac{\p u}{\p \nu} = \ldots = \frac{\p^{k-1} u}{\p\nu^{k-1}} = 0 \quad \mbox{on }\; \p\Omega;
\end{align}
or the Navier boundary conditions
\begin{eqnarray}\label{N}
u = \Delta u = \ldots = \D^{k-1} u = 0\quad \mbox{on }\; \p\Omega.
\end{eqnarray}
Here $\Omega\subset \mathbb{R}^{N}$ ($N > 2k$) is a bounded domain with
smooth boundary and $f$ is a $\mathrm{C}^{1}(\Omega\times
\mathbb{R})$ function that we will specify later.
Define
\begin{align}\label{quad}
  \Lambda_u (\phi):= \int_{\O}|D^k \phi|^2-f'(x,u)\phi^2\quad
  \mbox{for }\,\phi \in \Sigma_k
\end{align}
where
\begin{equation*}
  D^k = \left\{\begin{array}{ll}\nabla\D^{\frac{k-1}{2}}& \mbox{for $k$ odd};\\
\D^{\frac{k}{2}}& \mbox{for $k$ even}
\end{array}\right.
\end{equation*}
and
\begin{equation*}\label{}
  \Sigma_k:= \left\{\begin{array}{ll}
  H_0^k(\O) & \mbox{if we work with \eqref{D}};\\
\left\{\phi \in H^k(\O), \phi = \D \phi= . ..= \D^{[\frac{k-1}{2}]} \phi=0 \;\mbox{on }\; \p\O\right\} & \mbox{if we work with \eqref{N}.}
\end{array}\right.
\end{equation*}
The Morse index of a classical solution $u$ of $(P_k)$, denoted by $i(u)$
is defined as the maximal dimension of all subspaces of $\Sigma_k$ such
that $\Lambda_u(\phi) < 0$ in $\Sigma\setminus\{0\}.$ We say that $u$ is stable if its Morse index is equal to
zero. Our aim here is to get some explicit estimates of $u$ using its Morse index $i(u)$.

\medskip
We begin by presenting some assumptions on the nonlinearity $f$:
\begin{enumerate}
  \item [$(H_1)$] (superlinearity) There exists $\mu> 0$ such that
  $$f'(x,s)s^{2}\geq (1+\mu)f(x,s)s > 0,\quad \mbox{for} \; |s|>s_{0}, \;\; x\in \Omega.$$
  \item [$(H_2)$] (subcritical growth) There exists $0<\theta<1$ such that
  $$\frac{2N}{N-2k}F(x,s)\geq (1+\theta)f(x,s)s,\quad \mbox{for} \; |s|>s_{0}, \;\; x\in \Omega,$$
  where $F(x,s)=\displaystyle{\int_{0}^{t}}f(x,t)dt.$
  \item [$(H_3)$] There is a constant $C\geq 0$ such that
  $$|\nabla_x F(x,s)|\leq C(F(x,s)+1), \quad x\in \Omega.$$
\end{enumerate}
We say that $f$ satisfies $(H_i)$ in $\R_+$, if we have the assumption $(H_i)$ only for $s$ large enough.

\medskip
For the second order case, i.e.~$k = 1$, Bahri and Lions obtained in \cite{c} the estimates of solutions in $H_0^1(\O)$ for superlinear and subcritical growth $f$, by using the blow-up technique and the Morse index of the solutions. Motivated by \cite{c}, based on some local interior estimates and careful boundary estimates, Yang obtained in \cite{lec} the first explicit estimates of $L^{p}$ or
$L^{\infty}$ norm for solutions to $(P_1)$ via the Morse index. More precisely,
Yang proved that
\begin{taggedtheorem}{A} Let $f$ satisfy $(H_1)$-$(H_3)$, then there exist positive
constant $C$, $\alpha$ and $\beta$ such that any $u \in C^2(\O)\cap C(\overline{\O})$, solution of $(P_1)$ satisfies $$\int_{\Omega}|f(x,u)|^{p_{0}}dx \leq
C({i(u)}+1)^{\alpha}, \quad \|u\|_{L^{\infty}(\O)}\leq C (i(u)+1)^{\beta},$$
\mbox{where }
$$ p_{0} = 1+\frac{(1+\theta)(N-2)}{(1-\theta)N+2(1+\theta)},\quad \alpha = \left(\frac{3}{2}+\frac{3}{2+\mu}\right)\frac{(2+\mu)^{2}}{3\mu+\mu^{2}}$$
and $$\beta = \frac{2\alpha}{p_{0}N(2-p_{0})}\left[\frac{2}{N(2-p_{0})}-\frac{1}{p_{0}}\right]^{-1}.$$
\end{taggedtheorem}
Hajlaoui, Harrabi and Mtiri revised in \cite{HHF} the results of \cite{lec}, they
obtained similar $L^{\infty}$-estimate for solution to $(P_1)$. The proof in \cite{HHF} is more transparent, and it allows them to get a slightly better estimate for large dimension $N:$
\begin{taggedtheorem}{B} Let $f$ satisfy $(H_1)$-$(H_3)$, then there exist positive
constant $C$, $\alpha'$ and $\beta'$ such that any classical solution $u$ of $(P_1)$ satisfies
$$\displaystyle{\int_{\Omega}}|\nabla u|^2 dx\leq
C(i(u)+1)^{\alpha'},\quad \|u\|_{L^{\infty}}\leq C (i(u)+1)^{\beta'},$$
where
$$\alpha'=\frac{4}{\mu}+3 \quad \mbox{and}  \quad \beta'= \dfrac{3\mu+4}{3\mu
\theta}\times \frac{3N^2(1-\theta)+N(7\theta-4)-2\theta+12}{N(N-2)^2}.$$
\end{taggedtheorem}

In this paper, we will try to handle the polyharmonic equations. Let
$$
\left\{
\begin{array}{ll}
(-\Delta)^k u=f(x,u) & \mbox{in }\; \Omega;\\
\mbox{$u$ satisfies \eqref{D}}, & \mbox{if $k$ is odd};\\
\mbox{$u$ satisfies \eqref{N}}, & \mbox{if $k$ is even.}
\end{array}
\right.
\leqno{(E_k)}
$$
To simplify the presentation, we will concentrate on the cases $k=2$ and $k=3$, even we believe that the results should hold true for general $k \in \N$. 
We will provide some $L^{p}$ and $L^{\infty}$ estimates in polynomial growth
 function of the Morse index, for classical solutions of
$(E_2)$ and $(E_3)$, provided suitable conditions on $f$. As far as we know, it seems to be the first time 
that some explicit estimates are obtained for polyharmonic problems via the Morse index.

\medskip
As in \cite{HHF}, we shall employ a cut-off function with compact support to derive a variant of the Pohozaev identity. This device allows us to avoid the
 spherical integrals raised in \cite{lec}, which are very difficult to control, especially for the polyharmonic situations. Furthermore, under $(H_1)$-$(H_3)$, the local $L^2$-estimate of
 $\nabla u$ and $\D u$ via the Morse index seem also difficult to derive for the polyharmonic equation  than for $(P_1)$ the second order case. As in \cite{HHF}, we need to exhibit the explicit dependence on $i(u)$ (see Lemma \ref{l.2.3a} and lemma \ref{l.2.3b} below). The following are our main results.
\begin{thm}\label{main2}
If $u$ is a classical solution of $(E_2)$ with $f \geq 0$ satisfying $(H_1)$-$(H_3)$ in $\R_+$; or if $u$ is a classical solution of $(E_3)$ with $f$ satisfying $(H_1)$-$(H_3)$, then there exists a positive constant $C$ independent of $u$ such that
 $$\int_{\Omega}|f(x,u)|^{p_k} dx \leq C({i(u)}+1)^{\alpha_{k}}$$ where
$$p_{k}=\frac{2N}{N(1-\theta)+2k(1+\theta)} \quad \mbox{and} \quad \alpha_{k}=\frac{4k(\mu+1)}{\mu } \;\; \mbox{where $k =$ 2 or 3 respectively}.$$
\end{thm}
By setting up a standard boot-strap iteration, as $f$ has subcritical growth, we can proceed similarly as in the proof of Theorem 2.2 in \cite{lec} and claim that
\begin{thm}\label{main3}
If $u$ is a classical solution of $(E_2)$ with $f \geq 0$ satisfying $(H_1)$-$(H_3)$ in $\R_+$; or if $u$ is a classical solution of $(E_3)$ with $f$ satisfying $(H_1)$-$(H_3)$, then there exists a positive constant $C$ independent of $u$ such that (for $k = 2$ or $3$ respectively),
$$ \|u\|_{L^{\infty}(\O)}\leq C (i(u)+1)^{\beta_{k}}, \quad \mbox{where } \; \beta_{k} = \frac{2k\alpha_{k}}{p_{k}N(2-p_{k})}\left[\frac{2k}{N(2-p_{k})}-\frac{1}{p_{k}}\right]^{-1}, \;\;\alpha_{k}=\frac{4k(\mu+1)}{\mu },$$
and $p_k$ is defined in Theorem \ref{main2}.
\end{thm}

By assumptions $(H_1)$ and $(H_2)$ in $\R$ (resp. in $\R_+$), there exist two positive
constants $C_1$ and $C_2$ such that for $|s|$ large enough (resp. for $s$ large enough),
\begin{align}\label{B}
 \frac{(N-2k)(1+\theta)}{2N} f(x,s)s-C_1 \leq F(x,s)\leq \frac{1}{2+\mu}f(x,s)s + C_1,
\end{align}
  \begin{align}\label{a}
 f(x,s)s \geq C_{1}(|s|^{2+\mu} - 1)
  \end{align}
  and
 \begin{align}\label{A}
  |f(x,s)|\leq C_{2}\left(|s|^{\frac{N(1-\theta)+2k(1+\theta))}{(N-2k)(1+\theta)}}+1\right).
\end{align}
This paper is organized as follows : We give the proof of Theorem \ref{main2} for $k=2$ and $k=3$ respectively in sections 2 and 3. In the following, $C$ denotes always a generic positive constant independent of the solution $u$, even their
value could be changed from one line to another one.

\section{Proof for $k=2$}
\setcounter{equation}{0}
Here we will prove Theorem \ref{main2} for $k = 2$.

\subsection{Preliminaries}
Let $y \in \mathbb{R}^N$ and $R>0.$ Throughout the paper, we denote by $B_R(y)$
the open ball of center $y $ and radius $R$ and $\partial\Omega_R(y):=\partial\Omega \cap B_R(y)$. For $x\in
B_{R}(y)\cap\Omega,$ let $n:=x-y.$ We denote also
$$u_{j_i\cdots j_k} :=\frac{\partial^{k} u}{\partial x_{j_1}\partial x_{j_2}\cdots\partial x_{j_k}}.$$

\medskip
First of all, we have the following Pohozaev identity.
\begin{lem}\label{l.2.1}
Let $u$ be a classical solution to $(E_2)$. Let $\psi \in C_c^{2}(B_R(y))$. Then
\begin{align*}
&\; \frac{2N}{N-4}\int_\O F(x,u)\psi dx +
\frac{2}{N-4}\int_\O \nabla_{x} F(x,u)\cdot n\psi dx
-\int_\O (\Delta u)^{2}\psi dx \\
 =& \;-\frac{4}{N-4}\int_\O \Delta u \nabla^{2}u (\nabla \psi, n) dx + \frac{1}{N-4}\int_\O (\nabla\psi\cdot n)(\Delta u)^{2} dx\\
& \; -\frac{4}{N-4}\int_\O (\nabla u\cdot \nabla \psi)\Delta u dx
- \frac{2}{N-4}\int_\O (\nabla u\cdot n)\Delta u\Delta \psi dx \\
& \; -\frac{2}{N-4}\int_\O F(x,u)\nabla \psi\cdot n dx - \frac{2}{N-4}\int_{\partial\Omega_R(y)} \frac{\p \Delta u}{\p\nu}(\nabla u\cdot n)\psi d\sigma.
\end{align*}
\end{lem}
The proof is classical by multiplying the equation by $(n\cdot\nabla u)\psi$ and integration by parts, so we omit it.

\medskip
 To establish a global estimate, we will cover the domain
 $\O$ by small balls and obtain local estimates.  To be more precise, consider $$\Omega_{1,R}:=
\left\{x \in \O: \mbox{ dist}(x,\partial \O)> \frac{R}{2}\right\}\;\;\mbox{and}\;\;\Omega_{2,R}:=\left\{x \in \O: \mbox{ dist}(x,\partial \O)\leq \frac{R}{3}\right\},\quad \forall\;  R>0. $$
The main difficulty is the estimates of $u$ near the boundary, that is, in $\Omega_{2,R}$. We need to choose carefully the balls as in \cite{lec}. Indeed, we will take balls with center lying in
\begin{align}\label{Gamma}
\Gamma(R):= \left\{x\in \mathbb{R}^N\backslash \O:\;\mbox{dist}(x,\partial \O)=
 \frac{R}{20}\right\},
\end{align}
The domain $\O\backslash \Omega_{2,R}$ will be covered by balls with center lying in $\Omega_{1,R}.$
The following lemma is devoted to the control of the boundary term for $y \in \Gamma(R)$ in the above Pohozaev identity.
 \begin{lem}\label{l.2.2}
  There exists
$R_1 > 0$ depending on $\O$ such that if $f(x,u)\geq 0$ and $u$ is a classical solution of $(E_2)$, then for any $0< R \leq R_1$ and $y\in
\Gamma(R),$ there holds
\begin{align}
 \nonumber\int_{\partial\Omega_R(y)}\frac{\p \Delta u}{\p\nu}(\nabla u\cdot n)\psi d\sigma\geq0,
\end{align}
for any nonnegative function $\psi \in C_c^{2}(B_R(y)).$
 \end{lem}
 \textbf{Proof.} As in the proof of Lemma 2.2 of \cite{lec}, there exists $R_1 >0$ such that if $0<R \leq R_1$ and $y \in \Gamma(R)$ then $\nu \cdot n \leq 0$ for any $x \in \partial\Omega_R(y).$

 \smallskip
 As $f(x,u)\geq 0,$ the maximum principle implies that $-\D u\geq 0$ in $\O$ as $\D u = 0$ on $\p\O$, hence $u \geq 0$. Therefore $\frac{\p\D u}{\p\nu} \geq 0$ on $\p\O$ and $\nabla u \cdot n = (n\cdot v)\frac{\p u}{\p \nu}\geq 0$ on $\p\O$, so we obtain the claim. \qed

 \smallskip
Consequently, we get
 \begin{prop}\label{c.2.1}
There exists $R_0 > 0$ small who satisfies the following property: Let $u$
be a classical solution of $(E_2)$ with $f \geq 0$ verifying $(H_{1})$-$(H_{3})$ in $\R_+$. Then for any $0< R\leq R_0$, $y \in\Gamma(R)$ and $0 \leq \psi \in C_c^{4}(B_R(y))$, there holds
 \begin{align}
 \label{newl1}
 \begin{split}
& \int_\O f(x,u)u\psi dx + \int_\O (\D u)^2\psi dx \\
\leq& \; CR\|\nabla \psi\|_\infty \int_{A_{R,\psi}(y)}f(x,u) u dx + CR^{2}\int_{A_{R,\psi}(y)} |\nabla^2(u\nabla\psi)|^2 dx\\
& \;  + C\Big(1 + R\|\nabla\psi\|_\infty\Big)\| \D u\|_{L^{2}(A_{R,\psi}(y))}^{2}+C\Big(R^2 \|\nabla(\D\psi)\|_{\infty}^{2} + \|\Delta \psi\|_{\infty}^{2}\Big)\|u\|_{L^{2}(A_{R,\psi}(y))}^{2}\\
 & \; +CR^{2}\Big( \|\Delta \psi\|_{\infty}^2 +\frac{ 1}{R^{2}}\|\nabla\psi\|_\infty^2+\|\nabla^{2}\psi\|_{\infty}^{2}\Big)\| \nabla u\|_{L^{2}(A_{R,\psi}(y))}^{2}  + CR^{N},
\end{split}
\end{align}
where $$A_{R,\psi}(y) = B_{R}(y)\cap\Omega\cap\{\nabla \psi\neq0\}.$$
Moreover, for $y \in\Omega_{1,R}$, the above inequality holds true if we replace $R$ by $\frac{R}{2}.$
\end{prop}
\textbf{Proof.} Let $y \in \Gamma(R)$ with $R < R_1$ and $0 \leq \psi \in C_c^{4}(B_R(y))$. Using Lemmas \ref{l.2.1}--\ref{l.2.2}, $(H_{1})$-$(H_{3})$ and \eqref{B}, we obtain
\begin{align}
\label{xx.0.0}
\begin{split}
& (1+\theta)\int_\O f(x,u)u\psi dx - \int_\O (\D u)^2\psi dx \\
\leq & \;\frac{4}{N-4}
 \int_{A_{R,\psi}(y)}|\Delta u||\nabla^{2}u (\nabla \psi, n)| dx + \frac{1}{N-4}\int_{A_{R,\psi}(y)}(\Delta u)^{2}|\nabla\psi\cdot n| dx \\
& \; + \frac{4}{N-4}\int_{A_{R,\psi}(y)}|\Delta u||\nabla u\cdot \nabla \psi| dx
+\frac{2}{N-4}\int_{A_{R,\psi}(y)}|\Delta u||\nabla u\cdot n||\Delta \psi| dx \\
&\; + \frac{1}{(N-4)}\int_{A_{R,\psi}(y)}f(x,u)u|\nabla \psi\cdot n| dx + C R\int_{B_{R}(y)\cap\O}f(x,u)u \psi dx + CR^{N}.
\end{split}
\end{align}
A direct calculation implies that
 \begin{align*}
\nabla^2u(\nabla\psi, n) = \sum_{ij} u_{ij}\psi_in_j = \sum_{ij}(u\psi_i)_{ij}n_j - u\nabla(\D\psi)\cdot n - \D\psi(\nabla u \cdot n) - \nabla^2\psi(\nabla u, n).
\end{align*}
By the Cauchy-Schwarz inequality, there exists $C > 0$ such that
\begin{align}
\label{xxx.0.0}
\begin{split}
 \int_{A_{R,\psi}(y)}|\Delta u|| \nabla^{2}u (\nabla \psi, n)| dx
\leq& \; C\int_{A_{R,\psi}(y)}|\Delta u|^2 dx + CR^{2}\int_{A_{R,\psi}(y)} u^{2}|\nabla(\D\psi)|^{2} dx \\
 &+ C R^{2}\int_{A_{R,\psi}(y)} |\nabla^2(u\nabla\psi)|^2 dx\\
  & \; + CR^{2}\int_{A_{R,\psi}(y)}|\nabla u|^{2}\Big(\|\Delta \psi\|_\infty^2 + \|\nabla^2\psi\|_\infty^2\Big) dx.
\end{split}
\end{align}
On the other hand, recall that $u = \D u = 0$ on $\p\O$ and $\psi \in C_c^4(B_R(y))$, multiplying the equation $(E_2)$ by $u\psi$ and
integrating by parts, we get readily
\begin{align}\label{3.8}
\begin{split}
 \int_\O (\D u)^2\psi dx -\int_ \O f(x,u)u\psi dx &\leq C \int_{A_{R,\psi}(y)}|\D u|\Big[|\nabla u\cdot \nabla \psi| + |u||\Delta \psi| \Big] dx\\
 &\leq C\int_{A_{R,\psi}(y)}\Big[(\Delta u)^2 + |\nabla u\cdot\nabla\psi|^2 + (\Delta\psi)^2u^2\Big] dx.
 \end{split}
\end{align}
Remark that
\begin{align*}
\frac{\theta}{2} \int_\O (\D u)^2\psi dx + \frac{\theta}{2}\int_\O f(x,u)u\psi dx = & \; (1+\theta)\int_\O f(x,u)u\psi dx - \int_\O (\D u)^2\psi dx \\
& \; +\left(1 + \frac{\theta}{2}\right)\left[\int_{\O}(\D u)^2\psi dx - \int_{\O}f(x,u)u\psi dx\right].
\end{align*}
Fix $R_0 \in (0, R_1)$ such that $CR_0 < 1$. Combining \eqref{xx.0.0}-\eqref{3.8}, using again Cauchy-Schwarz inequality, there holds clearly \eqref{newl1}. The proof for $y \in \O_{1, R}$ is completely similar , so we omit it.\qed

\begin{rmk}
The key point in \eqref{newl1} is that the integral over the ball $B_R(y)\cap\O$ is now controlled by the integrals over 
the annuli type domain $A_{R, \psi}(y)$ when we work with suitable cut-off function $\psi$.
\end{rmk}

Let $R>0,$ $y \in \Omega_{1,R}\cup \Gamma(R)$, $0<a<b$. Denote $$A:=A_{a}^{b}=\{x\in \mathbb{R}^{N};\; a<|x-y|<b\}, \quad A_\rho:= A_{a+\rho}^{b-\rho} \;\; \mbox{for }\; 0 < \rho < \frac{b-a}{4}.\eqno{(*)}$$ We will use also the following classical estimates.
\begin{lem}\label{l.2.3a} There exists a constant
$C>0$ depending only on $N$ such that for any $u \in H^2(\O)\cap H_0^1(\O)$ and $0 < \rho < \min(1, \frac{b-a}{4}),$ we have
$$\|\nabla u\|_{L^2(A_\rho\cap \O)}^2 \leq C\left(\frac{1}{\rho^{2}}\| u\|_{L^2(A\cap \O)}^2+\|\D u\|_{L^2(A\cap \O)}^2\right).$$
\end{lem}
\begin{rmk}\label{rem2.1}
If $f$ satisfies $(H_1)$, using \eqref{a}, there holds
\begin{equation*}
\|u\|_{L^2(A\cap\O)}^2 \leq C\left(\int_{A\cap \O}f(x,u)u dx\right)^{\frac{2}{2+\mu}}+ C.
 \end{equation*}
\end{rmk}

\subsection{Estimation via Morse index}
Let $u$ be a solution to $(E_2)$ with $f \geq 0$ and finite Morse index $i(u)$. For $y \in \Gamma(R)\cup \O_{1, R}$, denote
\begin{align}
 \label{defaj}
 A_j=:A_{a_j}^{b_j}\quad \mbox{with}  \;\;a_j = \frac{2(j+i(u))}{4(i(u)+1)}R,
\;\;b_j = \dfrac{2(j+i(u))+1}{4(i(u)+1)}R, \quad 1\leq j \leq i(u)+1.
\end{align}
Fix a cut-off function $\Phi \in C^\infty(\R)$ such that $\Phi =1$ in $[0, 1]$ and supp$(\Phi) \subset ({-\frac{1}{2}}, \frac{3}{2})$.
Let $$\phi_{j}(x):= \Phi\left(\dfrac{4 (i(u)+1)|x-y|}{R} - 2j - 2i(u)\right).$$
Then for any $1\leq j \leq i(u)+1$, $\phi_j \in C_c^\infty(B_R(y))$,
\begin{align}
\label{newest1}
\phi_{j}(x)=1  \;\; \mbox{in }\; A_{j},\quad \|\nabla\phi_{j}\|_{\infty}\leq \frac{C}{R}(1+i(u))\quad \mbox{and} \quad
\|\Delta\phi_{j}\|_{\infty}\leq \frac{C}{R^2} (1+i(u))^{2}.
\end{align}

We prove the following lemma.
  \begin{lem}
\label{lemnew} Let $f$ satisfy $(H_1)$ and let $u$ be a smooth solution to $(E_2)$ with Morse index $i(u) < \infty$. Then for any $0 < R \leq R_0$, $y\in \Gamma(R)\cup\O_{1, R}$, there exists $j_{0}\in \{1,2,...,1+i(u)\}$ verifying
\begin{align}
\label{newest0}
\int_{A_{j_0}\cap\Omega}(\D
u)^{2} dx + \int_{A_{j_0}\cap\Omega}f(x,u)u dx \leq C
\left(\frac{1+i(u)}{R}\right)^{\frac{4\mu+ 8}{\mu}}.
\end{align}
\end{lem}
 \noindent{\bf Proof.} First, for $\epsilon \in (0, 1)$ and $\eta \in C^2(\R^N)$,
\begin{align*}
\int_{\O} [\Delta (u\eta)]^2 dx & = \int_{\O}\left(u \D \eta+2\nabla u\nabla \eta + \eta\D u\right)^{2} dx\\
& \leq \left(1+\frac{\epsilon}{2}\right)\int_{\O}(\D u)^2 {\eta}^2 dx + \frac{C}{\e}\int_{\O}u^2
  (\D\eta)^2 dx +\frac{C}{\epsilon} \int_{\O} |\nabla u|^2|\nabla \eta|^2 dx.
   \end{align*}
Using $\Delta(u^2)  = 2|\nabla u|^2  + 2u\D u$, there holds
\begin{align}
\label{newest5}
\int_{\O}|\nabla u|^2|\nabla \eta|^2 dx \leq \frac{1}{2}\int_{\O} u^2 \Delta(|\nabla \eta|^2) dx + \int_{\O}|u||\D
u| |\nabla \eta|^2 dx.
\end{align}
 Take $\eta = \zeta^m$ with $m > 2$, $\zeta \geq 0$ and apply Young's inequality, we get
  \begin{align}
  \label{newest4}
  \begin{split}
  \int_{\O}|u||\D
u||\nabla \zeta^{m}|^2 dx & = m^{2}\int_{\O} |u||\D
u||\nabla \zeta|^2\zeta^{2m-2} dx\\
& \leq \epsilon^{2} \int_{\O}(\D u)^2 \zeta^{2m} dx +
  C_{\epsilon, m}\int_{\O} u^2 |\nabla \zeta|^4\zeta^{2m-4} dx.
  \end{split}
\end{align}
Here $C_{\e, m}$ denotes a constant depending only on $\e$ and $m$. Therefore
\begin{align}
\label{newest2}
\int_{\O} [\Delta (u \zeta^{m})]^2 dx \leq(\epsilon+1)\int_{\O }(\D u)^2 \zeta^{2m} dx + C_{\epsilon, m}\int_{\O}u^2 \Big[ |\D\zeta|^2 +|\nabla \zeta|^4+|\Delta(|\nabla \zeta|^2)|\Big]\zeta^{2m-4} dx.
 \end{align}

Consider now the family of functions $\{u\phi_{j}^{m}\}_{1\leq j \leq i(u)+1}$, $m > 2$. With the definition of $\phi_j$, it's easy to see that different $\phi_j$ are supported by disjoint sets for different $j$, so they
are linearly independent as $u > 0$ in $\O$. Therefore, there must exist $j_{0}\in
\{1,2,...,1+i(u)\}$ such that $\Lambda_u(u \phi_{j_0}^m )\geq 0$ where $\Lambda$ is the quadratic form given by \eqref{quad}. Combining $\Lambda_u(u \phi_{j_0}^m )\geq 0$ with \eqref{newest1} and  \eqref{newest2}, we obtain
    \begin{align}\label{0.2}
\displaystyle{\int_{\Omega}}f'
(x,u)u^{2}\phi_{j_{0}}^{2m} dx -(1+\epsilon)\int_{\Omega}(\D u)^2
\phi^{2m}_{j_{0}} dx \leq
\frac{C_\e}{R^{4}}(1+i(u))^{4}\int_{\Omega}u^{2}\phi_{j_{0}}^{2m-4} dx.
\end{align}

Moreover, multiply the equation $(E_2)$ by $u\eta^2$ and
integrate by parts, we get, using \eqref{newest5}
\begin{align*}
& \; \int_{\O}\Big[(\D u)^2 \eta ^2 - f(x,u)u \eta^2\Big] dx\\
= & \; -4\int_{\O}\eta\D u \nabla u\cdot\nabla \eta dx -2\int_{\O}\eta u\D u\D \eta dx - 2\int_{\O}u\D
u|\nabla \eta|^2 dx\\
\leq & \; \e\int_\O (\D u)^2 \eta ^2 dx + C_\e\int_\O u^2(\Delta \eta)^2 dx + C_\e \int_\O |\nabla u|^2|\nabla\eta|^2 dx - 2\int_{\O}u\D
u|\nabla \eta|^2 dx\\
\leq & \; \e\int_\O (\D u)^2 \eta ^2 dx + C_\e\int_\O u^2\Big[(\Delta \eta)^2 + |\Delta(|\nabla \eta|^2)\Big]dx  + C_\e\int_{\O}|u\D
u||\nabla \eta|^2 dx.
\end{align*}
Take now $\eta = \phi_{j_0}^m$ with $m = 2 + \frac{2}{\mu} > 2$, there holds as for \eqref{newest4},
\begin{align*}
\int_\O |u\D u| |\nabla \eta|^2 dx \leq \e \int_\O (\Delta u)^2\phi_{j_0}^{2m} dx + C_\e\int_\O u^2 \phi_{j_0}^{2(m-2)} |\nabla \phi_{j_0}|^4 dx.
 \end{align*}
By \eqref{newest1}, we deduce then
\begin{align}\label{0.4}
(1-2\epsilon)\int_{\O}(\D
u)^2\phi_{j_0}^{2m}dx - \int_{\O}f(x,u)u\phi_{j_0}^{2m} dx \leq
\frac{C_\epsilon}{R^{4}}(1+i(u))^{4}\int_{\Omega}u^{2}\phi_{j_{0}}^{2m-4} dx.
\end{align}

Let $\e < \frac{1}{2}$, multiplying \eqref{0.4} by $\frac{1+2\e}{1-2\e}$, using \eqref{0.2} and $(H_1),$ we get
  \begin{align*}
\epsilon\int_{\O}(\D
u)^2\phi_{j_0}^{2m} dx + \left(\mu-\frac{4\epsilon}{1-2\epsilon}\right)\int_{\O}f(x,u)u\phi_{j_0}^{2m}dx\leq
\frac{C_\epsilon}{R^{4}}(1+i(u))^{4}\int_{\Omega}u^{2}\phi_{j_{0}}^{2m-4} dx +C_\e.
\end{align*}
Fix now $\epsilon <\min(2, \frac{\mu}{4 + 2\mu})$, there holds
\begin{align*}
\int_{\O}(\D u)^2\phi_{j_0}^{2m} dx + \int_{\O}f(x,u)u\phi_{j_0}^{2m} dx \leq
\frac{C}{R^{4}}(1+i(u))^{4}\int_{\Omega}u^{2}\phi_{j_{0}}^{2m-4} dx +C.
\end{align*}
Therefore, using \eqref{a} and $R \leq R_0$, for any $\epsilon'>0$,
  \begin{align*}
\int_{\Omega}(\D
u)^{2}\phi_{j_{0}}^{2m}dx + \int_{\Omega}uf(x,u)\phi_{j_{0}}^{2m} dx & \leq
C_{\epsilon'}\left(\frac{1+i(u)}{R}\right)^\frac{4\mu+8}{\mu} + C +
\e' \int_{\O}|u|^{\mu+2}\phi_{j_{0}}^{(m-2)(\mu+2)}dx \\
& \leq
C_{\epsilon'}\left(\frac{1+i(u)}{R}\right)^\frac{4\mu+8}{\mu} + C
\epsilon' \int_{\O}f(x,u)u\phi_{j_{0}}^{(m-2)(\mu+2)}dx \\
& =
C_{\epsilon'}\left(\frac{1+i(u)}{R}\right)^\frac{4\mu+8}{\mu}+ C
\epsilon' \int_{\O}f(x,u)u\phi_{j_{0}}^{2m} dx.
\end{align*}
For the last line, we used $(m-2)(\mu+2) = 2m$. Take $\epsilon' > 0$ small enough, the estimate \eqref{newest0} is proved. \qed

\subsection{Proof of Theorem \ref{main2} completed  }
Now, we are in position to prove Theorem \ref{main2} for $k=2$. Fix
$$R = R_0, \quad \rho:=\frac{R}{10(i(u)+1)}, \quad A_{j_0,\rho}:=A_{a_{j_0}+\rho}^{b_{j_0}-\rho} \subset A_{j_0} \;\mbox{be as in $(*)$}.$$
According to Lemmas \ref{l.2.3a}, \ref{lemnew} and Remark \ref{rem2.1}, there exists a positive constant $C$ independent of $y \in \Gamma(R)\cup \O_{1, R}$ such that
 \begin{align}\label{3.7}
 \|\D u\|^2_{L^2(A_{j_0,\rho}\cap\O)}+ \|\nabla u\|^2_{L^2(A_{j_0,\rho}\cap\O)}\leq C
(1+i(u))^\frac{4\mu+8}{\mu}.
\end{align}
Here, $a_{j_{0}}$ and $b_{j_{0}}$ are defined in \eqref{defaj} with $j_0$ given by Lemma \ref{lemnew}.

\medskip
Consider a cut-off function $\xi_{j_0} \in C_c^4(B_{b_{j_0}-\rho}(y))$ verifying
$\xi_{j_0}(x)\equiv 1$ in $B_{a_{j_0}+\rho}(y)$,
 with $$\|\nabla
\xi_{j_{0}}\|_{\infty}\leq \frac{C}{R}(1+i(u)), \quad \|\Delta \xi_{j_{0}}\|_{\infty}\leq
\frac{C}{R^{2}}(1+i(u))^2.$$
Applying Proposition \ref{c.2.1} with $\psi = \xi_{j_{0}}$, as $A_{R, \psi}(y) \subset A_{j_0, \rho}\cap \O$, we get
\begin{align}
\label{qssss}
\begin{split}
& \int_{\O}f(x,u)u\xi_{j_0} dx + \int_{\O}(\D u)^2\xi_{j_0} dx\\
 \leq & \;\ C(1+i(u))\int_{A_{j_0, \rho}\cap \Omega}\Big[(\D u)^2+f(x,u)u\Big] dx + C\int_{A_{j_0, \rho}\cap \Omega} |\nabla^2(u\nabla\xi_{j_{0}})|^2 dx\\
  &+ C(1+i(u))^6 \| u\|_{L^{2}(A_{j_0, \rho}\cap\O)}^{2} + C(1+i(u))^4\| \nabla u\|_{L^{2}(A_{j_0, \rho}\cap\O)}^{2}+ CR^{N}.
\end{split}
\end{align}
Since  $u\nabla\xi_{j_0}=0$ on $\partial\O,$ by standard elliptic theory, there exists $C_\O > 0$ depending only on $\Omega$ such that
 \begin{align}\label{0.255}
 \begin{split}
\int_{\Omega} |\nabla^2(u\nabla\xi_{j_{0}})|^2 dx & \leq
C_\O \int_{\Omega}|\D (u \nabla\xi_{j_0})|^2 dx\\ &= C_\O \int_{A_{j_0, \rho}\cap\O}|\D (u \nabla\xi_{j_0})|^2 dx\\
& \leq C\int_{A_{j_0, \rho}\cap\O} \Big[u^2 |\nabla(\D\xi_{j_0})|^2 + |\nabla u|^2|\nabla^{2} \xi_{j_0}|^2 + (\D u)^2| \nabla\xi_{j_0}|^2\Big] dx.
\end{split}
\end{align}
From \eqref{qssss}, \eqref{0.255}, we get the following inequality
\begin{align}
\label{PL}
\begin{split}
& \int_{\O}f(x,u)u\xi_{j_0} dx + \int_{\O}(\D u)^2\xi_{j_0} dx\\
 \leq & \;\ C(1+i(u))\int_{A_{j_0, \rho}\cap \Omega}\Big[(\D u)^2+f(x,u)u\Big] dx + C(1+i(u))^2\|\D u\|_{L^{2}(A_{j_0, \rho}\cap\O)}^{2}\\
  &+ C(1+i(u))^6 \| u\|_{L^{2}(A_{j_0, \rho}\cap\O)}^{2} + C(1+i(u))^4\| \nabla u\|_{L^{2}(A_{j_0, \rho}\cap\O)}^{2}+ CR^{N}.
\end{split}
\end{align}
On the other hand, using Remark \ref{rem2.1} and Lemma \ref{lemnew}, there holds
\begin{align}\label{okm}
 \|u\|_{L^2(A_{j_{0}}\cap\O)}^2 \leq C\left(\int_{A_{j_{0}}\cap \O}f(x,u)u dx\right)^{\frac{2}{2+\mu}} + C
 \leq C(1+i(u))^{\frac{8}{\mu}}.
 \end{align}
Combining \eqref{newest0}, \eqref{3.7}, \eqref{PL} and  \eqref{okm}, one obtains
\begin{align*}
 \int_{\O}f(x,u)u\xi_{j_0} dx + \int_{\O}(\D u)^2\xi_{j_0} dx
\leq C(1+i(u))^\frac{8\mu+8}{\mu}.
\end{align*}
As $\frac{R}{2}< a_{j_{0}}$ and $R = R_0$, we get then for any $y \in \Gamma(R)\cup \O_{1, R}$,
\begin{align*}
\int_{B_{\frac{R_0}{2}}(y)\cap \O} \Big[|\D u|^{2} + f(x,u)u\Big] dx \leq C
(1+i(u))^\frac{8\mu+8}{\mu}.
\end{align*}
By covering argument and \eqref{A}, we get finally
$$\int_{\O} f(x, u)^{p_2} dx \leq C\int_\O f(u)u dx + C \leq C
(1+i(u))^{\alpha_{2}},$$
where $p_{2}=\frac{2N}{N(1-\theta)+4(1+\theta)}$ and
$ \alpha_{2}=\frac{8(\mu+1)}{\mu }.$ So we are done.\qed

\section{Proof of Theorem \ref{main2} for $k=3$}
\setcounter{equation}{0}
In this section, we consider the equation $(E_3)$. We will proceed as for $(E_2)$ and keep the same notations, but we replace the Navier boundary conditions by the Dirichlet boundary conditions and we have no more the sign condition for $f$.

\subsection{Preliminaries}
We make some preparations here. For $\psi \in C^m$ for $m \geq 1$, to simplify the notation, we define
\begin{align*}
[\psi]_m(x) = \sum_{|\beta_1| + \ldots + |\beta_p| = m, |\beta_i|\geq 1} \prod_{i=1}^p |\p_{\beta_i}\psi(x)|
\end{align*}
and the semi-norms
\begin{align*}
|\psi|_{m, \infty} = \sum_{\alpha_1 + \ldots + \alpha_p = m, \alpha_i\geq 1} \prod_{i=1}^p \|\nabla^{\alpha_i}\psi\|_\infty, \forall \; m \geq 1.
\end{align*}
Obviously, for any $\psi \in C^m$, we have $\|[\psi]_m\|_\infty \leq C_m|\psi|_{m, \infty}$.

 \begin{lem}\label{l.2.7a}
Let $m \geq 3$. For any $\epsilon > 0$, there exists $C_{\e, m} >0$ such that for any $u \in H_0^3(\O)$ and $\zeta \in C^6(\overline\O)$, there holds
\begin{eqnarray}\label{PLa}
 \int_{\O}\Big[(\D u)^2|\nabla \zeta^{m}|^2 + |\nabla u|^2|\nabla^{2} \zeta^{m}|^2\Big] dx
\leq \e\int_{\O}|\nabla(\D u)|^2 \zeta^{2m} dx + C_\e\int_{\O} u^2 [\zeta]_6\zeta^{2m-6} dx.
\end{eqnarray}
\end{lem}
\noindent{\bf Proof.} Using the equality $\D(u^2) = 2u\D u + 2|\nabla u|^2$, we have
 \begin{align*}
  \int_{\O} |\nabla u|^{2} |\nabla \zeta|^4\zeta^{2m-4} dx & \leq \frac{1}{2}\int_{\O} u^2 \D\left(|\nabla \zeta|^4\zeta^{2m-4}\right) dx + \int_{\O}|u||\D u||\nabla \zeta|^4\zeta^{2m-4} dx.
\end{align*}
Applying Young’s inequality, we get, for any $\epsilon >0$
 \begin{align*}
 \int_{\O}|u\D u||\nabla \zeta|^4\zeta^{2m-4} dx
 \leq \e\int_{\O}(\D u)^2 |\nabla \zeta|^2\zeta^{2m-2} dx +
  C_\e\int_{\O} u^2 |\nabla \zeta|^6\zeta^{2m-6} dx.
\end{align*}
So we get
\begin{align}
\label{newest7}
\int_{\O} |\nabla u|^{2} |\nabla \zeta|^4\zeta^{2m-4} dx & \leq \e \int_{\O}(\D u)^2 |\nabla \zeta|^2\zeta^{2m-2} dx +
  C_\e\int_{\O} u^2 [\zeta]_6 \zeta^{2m-6} dx.
\end{align}

On the other hand, direct integrations by parts yield (recall that $u \in H_0^3(\O)$)
  \begin{align*}
 \int_{\O}(\D u)^2|\nabla\eta|^2 dx =&\; - \int_{\O}\nabla u\nabla(\D u)|\nabla \eta|^2 dx - 2\int_{\O} \D u \nabla^2\eta(\nabla \eta, \nabla u) dx\\
 =& \; - \int_{\O}\nabla u\nabla(\D u)|\nabla \eta|^2 dx + 2\int_{\O} u \nabla^2\eta(\nabla \eta, \nabla(\D u)) dx\\
 &\; +2\int_{\O}u\D u|\nabla^2\eta|^2 dx + 2\int_{\O}u\D u\nabla \eta \cdot \nabla(\D \eta) dx\\
   = &\; - \int_{\O}\nabla u\nabla(\D u)|\nabla \eta|^2 dx + 2\int_{\O} u \nabla^2\eta(\nabla \eta, \nabla(\D u)) dx\\
 &\; + \int_{\O}\Big[\D(u^2)- 2|\nabla u|^2\Big]|\nabla^2\eta|^2 dx + 2\int_{\O}u\D u\nabla \eta \cdot \nabla(\D \eta) dx.
\end{align*}
Hence
\begin{align}
\label{11.0.0a}
\begin{split}
  \int_{\O}\Big[(\D u)^2|\nabla \eta|^2 + 2|\nabla u|^2|\nabla^{2}\eta|^2\Big] dx = &\; -\int_{\O}\nabla u\nabla(\D u)|\nabla \eta|^2 dx + 2\int_{\O} u \nabla^2\eta(\nabla \eta, \nabla(\D u)) dx\\
 &\; + \int_{\O} u^2 \D(|\nabla^2\eta|^2) dx + 2\int_{\O}u\D u\nabla \eta \cdot \nabla(\D \eta) dx.
\end{split}
\end{align}
Consider $\eta = \zeta^m$. For any $\epsilon >0$, by Cauchy-Schwarz inequality, we have
\begin{align*}
& \; -\int_{\O}\nabla u\nabla(\D u)|\nabla \eta|^2 dx + 2\int_{\O} u \nabla^2\eta(\nabla \eta, \nabla(\D u)) dx\\
\leq & \; \e \int_{\O}  |\nabla(\D u)|^2\zeta^{2m} dx + C_\e \int_{\O} |\nabla u|^{2}|\nabla \zeta|^4\zeta^{2m-4} dx + C_{\epsilon}\int_{\O} u^2[\zeta]_6\zeta^{2m-6}dx
\end{align*}
and
\begin{align*}
2\int_{\O}u\D u\nabla \eta \cdot \nabla(\D \eta) dx
\leq \epsilon\int_{\O}|\D u|^2|\nabla \zeta^{m}|^2 dx + C_{\epsilon}\int_{\O} u^2[\zeta]_6\zeta^{2m-6}dx.
\end{align*}
Inserting the two above estimates in \eqref{11.0.0a}, one gets
\begin{align*}
& \; (1-\epsilon) \int_{\O}(\D u)^2|\nabla \zeta^{m}|^2 dx + \int_{\O}|\nabla u|^2|\nabla^{2} \zeta^{m}|^2 dx\\
\leq & \; \e \int_{\O}  |\nabla(\D u)|^2\zeta^{2m} dx + C_\e \int_{\O} |\nabla u|^{2}|\nabla \zeta|^4\zeta^{2m-4} dx + C_{\epsilon}\int_{\O} u^2[\zeta]_6\zeta^{2m-6}dx.
\end{align*}
Take another small enough $\e$ in \eqref{newest7}, there holds
\begin{align*}
& \; (1-2\epsilon) \int_{\O}(\D u)^2|\nabla \zeta^{m}|^2 dx +\int_{\O}|\nabla u|^2|\nabla^{2} \zeta^{m}|^2 dx \leq \e \int_{\O}  |\nabla(\D u)|^2\zeta^{2m} dx + C_{\epsilon}\int_{\O} u^2[\zeta]_6\zeta^{2m-6}dx.
\end{align*}
The proof is completed. \qed

\medskip
Using Lemma \ref{l.2.7a}, we obtain also
\begin{lem}\label{l.2.8a}
Let $m \geq 3$. For any $0 < \epsilon< 1,$ there exists
$C_\epsilon>0$ such that for any $u \in H_0^3(\O)$ and $\zeta \in C^6(\overline\O)$,
\begin{align*}
 \int_{\O} \Big[|\nabla u|^2(\D\zeta^m)^2 + |\nabla^{2} u|^2|\nabla \zeta^{m}|^2\Big] dx \leq \e\int_{\O} |\nabla(\D u)|^2\zeta^{2m} dx + C_{\epsilon}\int_{\O} u^2 [\zeta]_6\zeta^{2m-6} dx.
\end{align*}
\end{lem}

\noindent{\bf Proof.} From \eqref{newest5}, we obtain
\begin{align}
\label{aj.ja}
\begin{split}
  \int_{\O} |\nabla^{2} u|^2|\nabla \zeta^m|^2 dx \leq &\; \frac{1}{2} \int_{\O} |\nabla u|^2\D(|\nabla \zeta^m|^2) dx + m^2\int_{\O} |\nabla u\cdot\nabla(\D u)||\nabla \zeta|^2 \zeta^{2m-2}dx\\
\leq &\;\int_\O |\nabla u|^2\Big[C_\e|\nabla \zeta|^4\zeta^{2m-4} + \nabla \zeta^m\nabla(\D\zeta^m)\Big]dx + \int_\O |\nabla u|^2|\nabla^2\zeta^m|^2 dx \\
&\; + \e\int_\O |\nabla(\D u)|^2\zeta^{2m} dx.
\end{split}
\end{align}
Rewrite
$$C_\e|\nabla \zeta|^4\zeta^{2m-4} + \nabla \zeta^m\nabla(\D\zeta^m) = \zeta^{2m-4}\nabla \zeta \cdot\Psi$$
with a smooth function $\Psi$. In the spirit of \eqref{newest5}, we have
\begin{align}
\label{bvl}
\begin{split}
\int_\O |\nabla u|^2\zeta^{2m-4}\nabla \zeta \cdot\Psi dx & \leq \frac{1}{2} \int_{\O} u^2\D(\zeta^{2m-4}\nabla \zeta \cdot\Psi) dx + \int_{\O} |u||\D u|\zeta^{2m-4}\nabla \zeta \cdot\Psi dx\\
& \leq \int_\O u^2\Big[|\D(\zeta^{2m-4}\nabla \zeta \cdot\Psi)| + C_\e|\Psi|^2\zeta^{2m-6}\Big] dx +\e\int_{\O} (\Delta u)^2|\nabla\zeta|^2\zeta^{2m-2}dx\\
& \leq C_\e \int_\O u^2[\zeta]_6\zeta^{2m-6} dx + \e \int_{\O} (\Delta u)^2|\nabla\zeta^m|^2 dx.
\end{split}
\end{align}
Combining \eqref{PLa} and  \eqref{aj.ja}-\eqref{bvl}, there holds
 \begin{align*}
\int_{\O} |\nabla^{2} u|^2|\nabla \zeta^{m}|^2 dx \leq \e\int_{\O} |\nabla(\D u)|^2\zeta^{2m} dx + C_{\epsilon}\int_{\O} u^2 [\zeta]_6\zeta^{2m-6} dx.
\end{align*}
Furthermore, integrating by parts,
\begin{align*}
\int_{\O} |\nabla u|^2(\D\zeta^m)^2 dx = &\; -2\int_\O \nabla^2 u(\nabla u, \nabla\zeta^m)\D\zeta^m dx - \int_\O |\nabla u|^2\nabla(\D\zeta^m)\nabla\zeta^m dx\\
\leq & \; \frac{1}{2}\int_{\O} |\nabla u|^2(\D\zeta^m)^2 dx + C\int_\O|\nabla^2 u|^2|\nabla\zeta^m|^2 dx\\
& \; + \frac{1}{2}\int_\O u^2\D\Big[\nabla(\D\zeta^m)\nabla\zeta^m\Big] dx + \int_\O |u||\D u| \nabla(\D\zeta^m)\nabla\zeta^m dx.
\end{align*}
We deduce that
\begin{align*}
\int_{\O} |\nabla u|^2(\D\zeta^m)^2 dx \leq C\int_\O\Big[|\nabla^2 u|^2|\nabla\zeta^m|^2 + |\D u|^2|\nabla\zeta^m|^2\Big] dx + C\int_\O u^2[\zeta]_6\zeta^{2m-6} dx,
\end{align*}
so using the previous estimates, we are done.\qed

\medskip
Let $R>0,$ $y \in \Omega_{1,R}\cup \Gamma(R)$, $0<a<b$. Denote $A:=A_{a}^{b}$ and $A_\rho:= A_{a+\rho}^{b-\rho}$, similar to Lemma \ref{l.2.3a}, we have
\begin{lem}\label{l.2.3b} There exists a constant
$C>0$ depending only on $N$ such that for any $u \in H^3_0(\O)$ and $0 < \rho < \min(1, \frac{b-a}{4}),$ we have
$$\|\D u\|_{L^2(A_\rho\cap \O)}^2\leq C \left(\frac{1}{\rho^{4}}\|u\|_{L^2(A\cap\O)}^2+\|\nabla(\D u)\|_{L^2(A\cap \O)}^2\right).$$
\end{lem}

\subsection{Explicit estimate via Morse index}
  \begin{lem}\label{l.2.4a}
\label{lemnew3} Let $f$ satisfies $(H_1)$ and $u$ be a solution to $(E_3)$ with finite Morse index $i(u)$. Then for any $y\in \Gamma(R)\cup\O_{1, R}$ with $R > 0$, there exists $j_{0}\in \{1,2,...,1+i(u)\}$ such that
\begin{align}
 \nonumber \int_{A_{j_0}\cap\Omega}|\nabla(\D u)|^2 dx + \int_{A_{j_0}\cap\Omega}f(x,u)u dx \leq C
\left(\frac{1+i(u)}{R}\right)^\frac{6\mu+12}{\mu}.
\end{align}
\end{lem}
 \noindent{\bf Proof.}  Take $\eta \in C^6(\overline\O)$. By direct
 calculations, we get, as $u \in H_0^3(\O)$,
\begin{align*}
\int_{\O} [\nabla(\Delta (u\eta))]^2 dx = & \; \int_{\O}\left(\nabla(\Delta u)\eta +\D u \nabla \eta +2\nabla^{2} u \nabla \eta+\nabla u \D \eta +
  2\nabla u \nabla^{2} \eta+ u \nabla(\D \eta)\right)^{2}\\
 \leq &\; (1+\epsilon)\int_{\O} |\nabla(\Delta u)|^2\eta^2 dx\\
 & \; + C_\e
  \int_{\O} \Big[|\D u|^2|\nabla \eta|^2 + |\nabla^{2} u|^2|\nabla \eta|^2
  + |\nabla u|^2\left(|\nabla^{2} \eta|^2+|\D \eta|^2\right) + u^2 |\nabla(\Delta \eta)|^{2}\Big]dx.
\end{align*}
Using Lemmas \ref{l.2.7a}-\ref{l.2.8a}, let $\eta = \zeta^m$ with $m = 3 + \frac{6}{\mu} > 3$, we derive that
\begin{align*}
 \int_{\O} |\nabla(\Delta(u\zeta^{m}))|^2 dx
  \leq (1 + \e)\int_{\O} |\nabla(\Delta u)|^2\zeta^{2m} dx + C_\e \int_{\O} u^2[\zeta]_6\zeta^{2m-6} dx.
\end{align*}
As in section 2, we can easily check that $\{u\phi_{j}^{m}\}_{1\leq j \leq i(u) + 1}$
are linearly independent, so there exists $j_{0}\in
\{1,2,...,1+i(u)\}$ such that $\Lambda_u(u \phi_{j_0}^m )\geq 0$. The above estimate with $\zeta = \phi_{j_0}$ implies then
    \begin{align}\label{2.Da}
\displaystyle{\int_{\Omega}}f'
(x,u)u^{2}\phi_{j_{0}}^{2m}dx -(1+\epsilon)\int_{\Omega}|\nabla(\Delta u)|^2
\phi^{2m}_{j_{0}} dx \leq
\frac{C_\epsilon}{R^{6}}(1+i(u))^{6}\int_{\Omega}u^{2}\phi_{j_{0}}^{2m-6} dx.
\end{align}
Now, take $u \phi_{j_0}^{2m}$ as the test function for $(E_3)$, the integration by parts yields that
\begin{align*}
\int_{\O}|\nabla(\Delta u)|^2 \phi_{j_0}^{2m} dx - \int_{\O}f(x,u)u \phi_{j_0}^{2m} dx
= \int_\O \nabla(\D u)\cdot\Big[\nabla\left(\D(u\phi_{j_0}^{2m})\right) - \nabla(\D u)\phi_{j_0}^{2m}\Big] dx.
\end{align*}
Developing the right hand side, applying again  Lemmas \ref{l.2.7a}-\ref{l.2.8a}, we can conclude: For any $\epsilon>0$, there exists $C_\epsilon$
 such that
\begin{align}\label{2.La}
(1 - \epsilon)\int_{\O}|\nabla(\Delta u)|^2 \phi_{j_0}^{2m}dx-\int_{\O}f(x,u)u\phi_{j_0}^{2m}dx\leq
\frac{C_\epsilon}{R^{6}}(1+i(u))^{6}\int_{\Omega}u^{2}\phi_{j_{0}}^{2m-6} dx.
\end{align}
Multiplying \eqref{2.La} by $\frac{1+2\epsilon}{1-\epsilon}$
adding it with \eqref{2.Da}, we obtain  from $(H_1)$ that
  \begin{align*}
\epsilon\int_{\O}|\nabla(\Delta u)|^2 \phi_{j_0}^{2m} dx +\left(\mu-\frac{3\epsilon}{1-\epsilon}\right)\int_{\O}f(x,u)u\phi_{j_0}^{2m} dx \leq
\frac{C_\epsilon}{R^{6}}(1+i(u))^{6}\int_{\Omega}u^{2}\phi_{j_{0}}^{2m-6} dx +C.
\end{align*}
Fix $0 < \epsilon < \frac{\mu}{3+\mu}$, we get
\begin{align*}
 \int_{\O}|\nabla(\Delta u)|^2 \phi_{j_0}^{2m} dx +\int_{\O}f(x,u)u\phi_{j_0}^{2m} dx \leq
\frac{C}{R^{6}}(1+i(u))^{6}\int_{\Omega}u^{2}\phi_{j_{0}}^{2m-6} dx + C.
\end{align*}
By Young's inequality, for any $\epsilon' > 0$, there holds
  \begin{align*}
 \int_{\O}|\nabla(\Delta u)|^2 \phi_{j_0}^{2m} dx + \int_{\Omega}uf(x,u)\phi_{j_{0}}^{2m} dx & \leq
C_{\epsilon'}\left(\frac{1+i(u)}{R}\right)^\frac{6\mu+12}{\mu} +
\epsilon' \int_{\O}|u|^{\mu+2}\phi_{j_{0}}^{(m-3)(\mu+2)} dx\\
& \leq C_{\epsilon'}\left(\frac{1+i(u)}{R}\right)^\frac{6\mu+12}{\mu}+C
\epsilon' \int_{\O}f(x,u)u\phi_{j_{0}}^{2m} dx.
\end{align*}
We used \eqref{a} and $(m-3)(2+\mu) = 2m$ for the last line. Take $\epsilon'$ small enough, the claim follows.\qed

\subsection{Proof of Theorem \ref{main2} for $k = 3$}
We show firstly the Pohozaev identity  for $(E_3)$.
\begin{lem}\label{l.2.1a}
Let $u$ be solution to $(E_3)$. Let $\psi \in C_c^{4}(B_R(y))$. Then
\begin{align*}
& N\int_\O F(x,u)\psi dx +
 \int_\O \nabla_{x} F(x,u)\cdot n\psi dx
- \frac{N-6}{2}\int_\O |\nabla(\Delta u)|^{2}\psi dx \\
 =&\; \frac{1}{2}\int_\O |\nabla(\Delta u)|^{2}(\nabla\psi\cdot n) dx  - \int_\O F(x,u)\nabla \psi\cdot n dx \\
&
\;- \int_\O \D\psi\nabla(\Delta u)\nabla(n\cdot\nabla u) dx - 2\int_\O \nabla(\Delta u)\nabla\Big[\nabla^2u(n, \nabla\psi) + \nabla u\nabla\psi\Big] dx\\
&+ \int_{\partial\Omega_R(y)}\frac{\p\Delta u}{\p \nu}\left(\nabla(\Delta u)\cdot n\right)\psi d\sigma -\frac{1}{2}\int_{\partial\Omega_R(y)}|\nabla(\Delta u)|^{2}(\nu\cdot n)\psi d\sigma.
\end{align*}
\end{lem}

For the boundary terms, we have
\begin{lem}\label{l.2.2a}
 There exists $R_1 > 0$ depending only on $\O$ such that for any $u$ smooth function in $H_0^3(\O)$, any $0<R<R_1$, $y\in
\Gamma(R)$ and any nonnegative function $\psi$, there holds
\begin{align}
 \nonumber \int_{\partial\Omega_{R}(y)}\frac{\p\Delta u}{\p \nu}(\nabla(\Delta u)\cdot n)\psi d\sigma -\frac{1}{2}\int_{\partial\Omega_{R}(y)}|\nabla(\Delta u)|^{2} \nu\cdot n\psi d\sigma \leq 0.
\end{align}
 \end{lem}

\noindent
\textbf{Proof.} Take $R_1 >0$ such that $v\cdot n \leq 0$ on $\p\O_R(y)$ for any $0<R\leq R_1$ and $y \in \Gamma(R)$. As $u \in H_0^3(\O)$, we know that $\nabla(\D u)$ is parallel to $\nu$ on $\p\O$, in other words $\nabla(\D u)(x) = \lambda(x)\nu(x)$ on $\p\O$. Therefore
 \begin{align*}
\frac{\p\Delta u}{\p \nu}(\nabla(\Delta u)\cdot n) - \frac{1}{2}(\nu\cdot n)|\nabla(\Delta u)|^2 = \frac{\lambda^2}{2}(\nu\cdot n)\leq 0,\quad \forall\; x \in \p\O_R(y).
 \end{align*}
So we are done. \qed

\medskip
Similar to Proposition \ref{c.2.1}, we can claim
 \begin{prop}\label{c.2.1a}
There exists $R_0 > 0$ small who satisfies the following property: Let $u$
be a classical solution of $(E_3)$ with $f$ verifying $(H_{1})$-$(H_{3})$. Then for any $0< R\leq R_0$, $y \in\Gamma(R)$ and $\zeta \in C_c^6(B_R(y))$ verifying $0 \leq \zeta \leq 1$ and $\psi = \zeta^{2m}$ with $m \geq 3$, there holds
\begin{align}
\label{xx.0.ve}
 \begin{split}
 &\; \int_\O f(x,u)u\psi dx + \int_\O |\nabla(\Delta u)|^2\psi dx\\
\leq& \; CR\|\nabla \zeta\|_\infty \int_{A_{R,\psi}(y)}f(x,u) u dx + C\Big(1 + R\|\nabla \zeta\|_\infty + R^2|\zeta|_{2, \infty}\Big)\|\nabla(\D u)\|_{L^{2}(A_{R,\psi}(y))}^{2}\\
 & \; + CR^{2}|\zeta|_{6, \infty}\|\nabla u\|_{L^{2}(A_{R,\psi}(y))}^{2} + C\Big(|\zeta|_{6, \infty} + R^2|\zeta|_{8, \infty}\Big) \|u\|_{L^{2}(A_{R,\psi}(y))}^{2}.
\end{split}
\end{align}
\end{prop}
\textbf{Proof.} Using Lemmas \ref{l.2.1a}- \ref{l.2.2a},  $(H_{1})$--$(H_{3})$  and  by $\eqref{B}$, we obtain
\begin{align}
\label{xx.0.0Ava}
\begin{split}
 &\; \frac{N-6}{2}\left[(1+\theta)\int_{\O}f(x,u)u\psi dx - \int_{\O}|\nabla(\Delta u)|^2\psi dx \right]\\
\leq&\; CR\|\nabla\psi\|_\infty \int_{A_{R,\psi}(y)}|\nabla(\Delta u)|^2 dx + CR\int_\O f(x, u)u\psi dx + CR\|\nabla\psi\|_\infty\int_{A_{R,\psi}(y)} f(x, u) u dx \\
&\; + \int_{A_{R,\psi}(y)} \Big|\nabla(\Delta u)\nabla\Big[\nabla^2u(n, \nabla\psi) + \nabla u\nabla\psi\Big]\Big| dx + \int_\O \Big|\D\psi\nabla(\Delta u)\nabla(n\cdot\nabla u)\Big| dx + CR^{N}.
\end{split}
\end{align}

 \noindent
 We will use also the following lemma.
\begin{lem}\label{l.2.7aa}
For any  $R < 1$, $\psi = \zeta^{2m}$ with $\zeta \in C_c^6(B_R(y))$ in Proposition \ref{c.2.1a}, there exists a positive constant
$C$ such that
\begin{align}
\label{xxx.013.0Ba}
\begin{split}
&\; \int_{A_{R,\psi}(y)} \Big|\nabla(\Delta u)\nabla\Big[\nabla^2u(n, \nabla\psi) + \nabla u\nabla\psi\Big]\Big| dx + \int_\O \Big|\D\psi\nabla(\Delta u)\nabla(n\cdot\nabla u)\Big| dx\\
 \leq& \;C \int_{A_{R,\psi}(y)}|\nabla(\Delta u)|^{2}dx+ CR^{2}\int_{A_{R,\psi}(y)} |\nabla(\D u)|^{2}[\zeta]_2 dx \\
 & \; + CR^{2}\int_{A_{R,\psi}(y)} |\nabla u|^2[\zeta]_6 dx + \int_{A_{R,\psi}(y)} u^2\left([\zeta]_6 + R^2[\zeta]_8\right) dx.
\end{split}
\end{align}
\end{lem}
\noindent{\bf Proof.} Indeed, in $B_R(y)\cap\O$,
 \begin{align*}
\Big|\nabla(\Delta u)\nabla\Big[\nabla^2u(n, \nabla\psi) + \nabla u\nabla\psi\Big]\Big|
 \leq & \; CR|\nabla(\D u)|\Big(|\nabla^3(u\nabla\psi)| + |\nabla^2u||\nabla^2\psi| + |\nabla u||\nabla^3\psi| + |u||\nabla^4\psi|\Big)\\
  &\; + C|\nabla(\D u)|\Big(|\nabla^2 u||\nabla\psi| + |\nabla u||\nabla^2\psi|\Big).
\end{align*}
We get then
\begin{align*}
& \int_{A_{R,\psi}(y)}\Big|\nabla(\Delta u)\nabla\Big[\nabla^2u(n, \nabla\psi) + \nabla u\nabla\psi\Big]\Big| dx \\
 \leq & \;C\int_{A_{R,\psi}(y)}|\nabla(\D u)|^2 dx + CR^2\int_{A_{R,\psi}(y)}|\nabla^3(u\nabla\psi)|^2 dx + CR^2\int_{A_{R,\psi}(y)}|\nabla^2u|^2|\nabla^2\psi|^2 dx \\
 & \; + CR^2\int_{A_{R,\psi}(y)}|\nabla u|^2|\nabla^3\psi|^2 dx + CR^2\int_{A_{R,\psi}(y)}u^2|\nabla^4\psi|^2 dx \\
  &\; + C\int_{A_{R,\psi}(y)}|\nabla^2 u|^2|\nabla\psi|^2 dx + C\int_{A_{R,\psi}(y)}|\nabla u|^2|\nabla^2\psi|^2 dx.
\end{align*}

First, using Lemmas \ref{l.2.7a}-\ref{l.2.8a} on $A_{R, \psi}(y)\cap \O$, the last two terms can be upper bounded by
\begin{align*}
C\int_{A_{R,\psi}(y)}|\nabla(\D u)|^2 dx + C\int_{A_{R,\psi}(y)} u^2[\zeta]_6 dx.
\end{align*}
Moreover, as $u\nabla\psi \in H_0^3(\O)$, there exists $C>0$ depending only on $\O$ such that
\begin{align*}
\int_{A_{R,\psi}(y)}|\nabla^3(u\nabla\psi)|^2 dx = \int_\O|\nabla^3(u\nabla\psi)|^2 dx \leq C\int_\O|\nabla\D(u\nabla\psi)|^2 dx = C\int_{{A_{R,\psi}(y)}} |\nabla\D(u\nabla\psi)|^2 dx.
\end{align*}
Remark that (as $\psi = \zeta^{2m}$)
\begin{align*}
|\nabla\D(u\nabla\psi)|^2 & \leq C\Big(|\nabla(\D u)|^2|\nabla\psi)|^2 + |\nabla^2 u|^2|\nabla^2\psi|^2 + |\nabla u|^2|\nabla^3\psi|^2 + u^2|\nabla^4\psi|^2\Big)\\
& \leq C\Big(|\nabla(\D u)|^2[\zeta]_2 + |\nabla u|^2[\zeta]_6 + u^2[\zeta]_8\Big) + C|\nabla^2 u|^2|\nabla^2\psi|^2.
\end{align*}
Using the equality $2|\nabla^{2}u|^2 = \D(|\nabla u|^2)- 2\nabla u\cdot\nabla(\D u)$, we obtain
\begin{align}
\label{newest9}
\begin{split}
  \int_{A_{R,\psi}(y)} |\nabla^{2} u|^2|\nabla^{2} \psi|^2 dx \leq& \; \frac{1}{2} \int_{A_{R,\psi}(y)} |\nabla u|^2\D(|\nabla^{2} \psi|^2) dx + \int_{A_{R,\psi}(y)} |\nabla u\cdot\nabla(\D u)||\nabla^{2} \psi|^2 dx\\
\leq & \; \frac{1}{2} \int_{A_{R,\psi}(y)} |\nabla u|^2|\D(|\nabla^{2} \psi|^2)| dx + C\int_{A_{R,\psi}(y)}|\nabla(\Delta u)|^{2}|\nabla^{2} \psi| dx \\
 &\; +C\int_{A_{R,\psi}(y)}|\nabla u|^{2}|\nabla^{2} \psi|^3 dx\\
 \leq & \; \int_{A_{R,\psi}(y)} |\nabla u|^2[\zeta]_6 dx + C\int_{A_{R,\psi}(y)}|\nabla(\Delta u)|^{2}[\zeta]_2 dx.
 \end{split}
\end{align}
Combining all these inequalities, we obtain the estimate for the first left term in \eqref{xxx.013.0Ba}.

\medskip
On the other hand,
\begin{align*}
& \;\int_{A_{R,\psi}(y)}\Big|\D\psi\nabla(\Delta u)\nabla(n\cdot\nabla u)\Big| dx\\
 \leq & \; \int_{A_{R,\psi}(y)}|\nabla(\D u)|\Big[R|\nabla^2 u||\Delta \psi| + |\nabla u||\Delta\psi|\Big] dx \\
\leq & \; \int_{A_{R,\psi}(y)}|\nabla(\D u)|^2 dx + C\int_{A_{R,\psi}(y)}\Big[R^2|\nabla^2 u|^2|\nabla^2 \psi|^2 + |\nabla u|^2(\D\psi)^2\Big] dx.
\end{align*}
Applying \eqref{newest9} and Lemma \ref{l.2.8a}, the proof is completed.\qed

\medskip
Coming back to the proof of \eqref{xx.0.ve} . Take $u \zeta^{2m}$ as the test function for $(E_3)$, using Lemmas \ref{l.2.7a}-\ref{l.2.8a}, for any $\epsilon>0$ there exists $C_\epsilon$ such that
\begin{align}\label{xx.0.0Aa}
\int_{\O}|\nabla(\Delta u)|^2 \zeta^{2m} dx -\int_{\O}f(x,u)u\zeta^{2m} dx \leq \e\int_{\O}|\nabla(\Delta u)|^2 \zeta^{2m} dx + C_\e\int_{\Omega}u^{2}[\zeta]_6\zeta^{2m-6} dx.
\end{align}
Remark that
\begin{align*}
\frac{\theta}{2} \int_\O |\nabla(\Delta u)|^2 \psi dx + \frac{\theta}{2}\int_\O f(x,u)u\psi dx = & \; (1+\theta)\int_\O f(x,u)u\psi dx - \int_\O |\nabla(\Delta u)|^2\psi dx \\
& \; +\left(1 + \frac{\theta}{2}\right)\left[\int_{\O}|\nabla(\Delta u)|^2\psi dx - \int_{\O}f(x,u)u\psi dx\right].
\end{align*}
Combining \eqref{xx.0.0Ava}-\eqref{xxx.013.0Ba} and \eqref{xx.0.0Aa}, for $\e, R > 0$ small enough, we have \eqref{xx.0.ve}. \qed

\medskip \noindent
{\bf Proof of Theorem \ref{main2} for $k=3$ completed}.

 Now, we are in position to prove Theorem \ref{main2} for $k=3$ . Fix $$R = R_0, \quad m = 3+\frac{6}{\mu}, \quad \rho:=\frac{R}{10(i(u)+1)}, \quad A_{j_0,\rho}:=A_{a_{j_0}+\rho}^{b_{j_0}-\rho} \subset A_{j_0} \;\mbox{be as in $(*)$}.$$
Using Remark \ref{rem2.1} and  lemma \ref{l.2.4a}, there holds
\begin{eqnarray}\label{PPLj}
 \|u\|_{L^2(A_{j_{0}}\cap\O)}^2 \leq C\left(\int_{A_{j_{0}}\cap \O}f(x,u)u\right)^{\frac{2}{2+\mu}} + C
 \leq C(1+i(u))^{\frac{12}{\mu}}.
 \end{eqnarray}
According to Lemmas \ref{l.2.3a}, \ref{l.2.3b}, \ref{l.2.4a} and \eqref{PPLj}, there exists a positive constant $C$ independent of $y \in \Gamma(R)\cup \O_{1, R}$ such that
 \begin{align}\label{3.7b}
 \|\nabla(\D u)\|^2_{L^2(A_{j_0,\rho}\cap\O)}+ \|\nabla u\|^2_{L^2(A_{j_0,\rho}\cap\O)}\leq C
(1+i(u))^\frac{6\mu+12}{\mu}.
\end{align}
Combining \eqref{xx.0.ve}, \eqref{PPLj} and  \eqref{3.7b}, one obtains
\begin{align*}
 \int_{\O}f(x,u)u\xi_{j_0} dx + \int_{\O}(\D u)^2\xi_{j_0} dx
\leq C(1+i(u))^\frac{12\mu+12}{\mu}.
\end{align*}
As $\frac{R}{2}< a_{j_{0}}$ and $R = R_0$, we get then for any $y \in \Gamma(R)\cup \O_{1, R}$,
\begin{align*}
\int_{B_{\frac{R_0}{2}}(y)\cap \O} \Big[|\D u|^{2} + f(x,u)u\Big] dx \leq C
(1+i(u))^\frac{12\mu+12}{\mu}.
\end{align*}
The proof is completed by the covering argument.\qed

\bigskip\noindent

\textbf{Acknowledgment:} A.H. and F.M. would like to express their deepest gratitude to our Research Laboratory \emph{LR11ES53 Algebra, Geometry and Spectral Theory (AGST) Sfax University}, for providing us with an excellent atmosphere
for doing this work.

\end{document}